\documentclass[9pt,conference]{IEEEtran}
\usepackage{cite}
\usepackage{amsmath,amssymb,amsfonts}
\usepackage{algorithmic}
\usepackage{graphicx}
\usepackage{textcomp}
\usepackage{xcolor}
\usepackage{url}
\usepackage{breakurl}
\usepackage{upgreek}
\def\BibTeX{{\rm B\kern-.05em{\sc i\kern-.025em b}\kern-.08em
    T\kern-.1667em\lower.7ex\hbox{E}\kern-.125emX}}

\usepackage{bm}               
\usepackage{bbm}

\usepackage{amsthm}
\usepackage{amscd}
\usepackage{amsmath}
\usepackage{amssymb}

\usepackage{mathtools}
\usepackage{mathrsfs}

\usepackage{tikz}
\usepackage{graphicx}
\usepackage{epstopdf}

\usepackage{caption}
\usepackage{subcaption}

\usepackage[capitalise]{cleveref}

\usepackage{times}
\usepackage{xspace}
\usepackage{enumitem}
\usepackage{tabularx}

\usepackage{algorithm}
\usepackage{algorithmic}

\usepackage{comment}

\newcommand{\tk}[1]{\textcolor{red}{\text{[TK: }#1\text{]}}} 

\newcommand{\vm}[1]{\textcolor{orange}{\text{[VM: }#1\text{]}}}

\newtheorem{theorem}{Theorem}[section]

\newtheorem{lemma}[theorem]{Lemma}
\newtheorem{corollary}[theorem]{Corollary}

\newtheorem{definition}[theorem]{Definition}

\newtheorem{problem}[theorem]{Problem}

\newtheorem{remark}[theorem]{Remark}

\newcommand{\ie}{\textit{i.e.}}

\newcommand{\e}{\mathrm{e}}             

\newcommand{\suml}{\sum\nolimits}

\newcommand{\half}{\dfrac{1}{2}} 

\newcommand{\beq}[1]{\begin{equation}\label{eq:#1}}
\newcommand{\eeq}{\end{equation}}

\newcommand{\defeq}{\coloneqq}                      
\newcommand{\eqdef}{\eqqcolon}                      

\newcommand{\inn}{\!\in\!}
\newcommand{\+}{\!+\!}
\renewcommand{\-}{\!-\!}
\renewcommand{\=}{\!=\!}

\renewcommand{\>}{\!>\!}
\newcommand{\<}{\!<\!}
\renewcommand{\leqq}{\!\leq\!}
\renewcommand{\geqq}{\!\geq\!}

\newcommand{\pr}[1]{\left({#1}\right)}          
\newcommand{\br}[1]{\left[{#1}\right]}          
\newcommand{\cl}[1]{\left\{{#1}\right\}}        

\newcommand{\abs}[1]{\vert{#1}\vert}                    

\newcommand{\norm}[2][\text{}]{\Vert{#2}\Vert_{#1}}     
\newcommand{\Norm}[2][\text{}]{\left\Vert{#2}\right\Vert_{#1}} 

\newcommand{\clf}[1]{\mathcal{#1}} 

\newcommand{\N}{\mathbb{N}}     
\newcommand{\Z}{\mathbb{Z}}     
\newcommand{\R}{\mathbb{R}}     

\newcommand{\Sym}{\mathbb{S}}   

\renewcommand{\j}{\jmath} 

\newcommand{\ejw}{\e^{\j\omega}}

\DeclareMathOperator*{\argmin}{arg\,min}

\newcommand{\subjto}{\,\mathrm{s.t.}\,}

\newcommand{\sub}{\mathscr{S}}

\newcommand{\tr}{\operatorname{tr}}         
\newcommand{\Tr}{\operatorname{Tr}}         

\newcommand{\tp}{\intercal}                 

\newcommand{\inv}{\mathrm{\-1}}             

\newcommand{\psdleq}{\preccurlyeq} 
\newcommand{\psdgeq}{\succcurlyeq} 
\newcommand{\psdl}{\prec}          
\newcommand{\psdg}{\succ}          

\newcommand{\E}{\operatorname{\mathbb{E}}} 

\newcommand{\Lp}{\mathscr{L}} 
\newcommand{\Hp}{\mathscr{H}} 

\newcommand{\F}{\mathcal{F}} 
\newcommand{\G}{\mathcal{G}} 
\newcommand{\K}{\mathcal{K}} 
\newcommand{\T}{\mathcal{T}} 
\newcommand{\NN}{\mathcal{N}}
\newcommand{\M}{\mathcal{M}} 
\renewcommand{\L}{\mathcal{L}} 
\newcommand{\RR}{\mathcal{R}}  

\newcommand{\I}{\mathcal{I}} 
\newcommand{\Ss}{\mathcal{S}}

\renewcommand{\u}{\mathbf{u}} 
\newcommand{\w}{\mathbf{w}} 
\newcommand{\x}{\mathbf{x}} 
\newcommand{\s}{\mathbf{s}} 

\newcommand{\QQ}{\clf{Q}}
\newcommand{\PP}{\clf{P}}
\newcommand{\Htwo}{\mathit{H}_2}
\newcommand{\Hinf}{\mathit{H}_\infty}
\newcommand{\RO}{\mathrm{RO}}

\newcommand{\FixAlg}{\texttt{} }

\newcommand{\causal}{\mathscr{K}}

\newcommand{\trig}{\mathscr{T}}

\newcommand{\HS}{2} 

\newcommand{\Policy}{\mathscr{K}}
\newcommand{\cost}{\operatorname{cost}}

\newcommand{\op}{\infty} 

\newcommand{\nc}{\textrm{nc}} 

\newcommand{\ght}{\gamma_{2}}
\newcommand{\ghi}{\gamma_{\infty}}

\newcommand{\Ktwo}{\K_{2}}
\newcommand{\Kinf}{\K_{\infty}}

\begin{document}

\title{Optimal Infinite-Horizon Mixed $\mathit{H}_2/\mathit{H}_\infty$ Control}

\author{
    \IEEEauthorblockN{
        Vikrant Malik\IEEEauthorrefmark{1}\IEEEauthorrefmark{2}\hspace{1cm} 
        Taylan Kargin \IEEEauthorrefmark{1}\IEEEauthorrefmark{2} \hspace{1cm}
        Joudi Hajar\IEEEauthorrefmark{1} \hspace{1cm}
        Babak Hassibi \IEEEauthorrefmark{1}
    } \\
    \IEEEauthorblockA{
        \IEEEauthorrefmark{1}\textit{Caltech, Pasadena, CA 91125, USA,}\hspace{1cm}
        \IEEEauthorrefmark{2}\textit{Equal contribution}
    }
}

\maketitle

\begin{abstract}
We study the problem of mixed $\Htwo/\Hinf$ control in the infinite-horizon setting. We identify the optimal causal controller that minimizes the $\Htwo$-cost of the closed-loop system subject to an $\Hinf$ constraint. Megretski \cite{megretski1994order} proved that the optimal mixed $\Htwo/\Hinf$ controller is non-rational whenever the constraint is active, without giving an explicit construction of the controller. In this work, we provide the first exact closed-form solution to the infinite-horizon mixed $\Htwo/\Hinf$ control in the frequency domain. While the optimal controller is non-rational, our formulation provides a finite-dimensional parameterization of the optimal controller. Leveraging this fact, we introduce an efficient iterative algorithm that finds the optimal causal controller in the frequency domain. We show that this algorithm is convergent when the system is scalar and present numerical evidence for exponential convergence of the proposed algorithm. Finally, we show how to find the best (in $H_\infty$ norm) fixed-order rational approximations of the optimal mixed $\Htwo/\Hinf$ controller, and study its performance. 

\end{abstract}

\section{Introduction}
Performance and robustness are the two most desired characteristics of a controller. Tackling uncertainty plays a pivotal role in the realm of control, especially within control systems that are frequently exposed to a range of uncertainties including external disruptions, inaccuracies in measurements, deviations in models, and time-varying system dynamics. Ignoring such uncertainties during the development of control policies can significantly undermine performance. The traditional $\Htwo$ \cite{kalman1960new} and $\Hinf$  \cite{doyle1988state} control are the two main approaches to address the dichotomy of robustness and performance. While $\Htwo$ control aims to achieve optimal average case performance for stochastic white disturbances, $\Hinf$ control guarantees robustness to worst-case deterministic disturbances with bounded energy. They both result in rational controllers which can be efficiently synthesized by solving a set of algebraic Riccati equations.


When it comes to designing a controller that is both robust to worst-case disturbances and optimal in the average case performance, a natural idea is to find the best $\Htwo$-optimal controller among a family of suboptimal $\Hinf$ controllers. Also known as the mixed $\Htwo/\Hinf$ control, these controllers provide a trade-off between the performance of $\Htwo$ and the robustness of the $\Hinf$ controllers. Although both $\Htwo$ and central $\Hinf$ controllers are rational and admit finite-order state-space realizations \cite{doyle1988state}, Megretski \cite{megretski1994order} showed that the mixed $\Htwo/\Hinf$ controller has infinite order whenever the $\Hinf$-norm constraint is not redundant.

\subsection{Prior Works}
The mixed $\Htwo/\Hinf$ control problem and its several variations have been extensively studied in the past. \cite{mustafa1988controllers} is one of the earliest works establishing a connection between the maximum entropy $\Hinf$ solution and the $\Htwo$ performance. \cite{bernstein1989lqg} proposed the first conceptualization and formulation of mixed $\Htwo$ control with $\Hinf$ constraints, and obtained a fixed order dynamic output feedback controller for an auxiliary objective upper bounding the $\Htwo$-norm. \cite{zhou1990mixed, khargonekar1991mixed} approached this problem similarly by considering auxiliary systems or objectives to bound the desired mixed $\Htwo/\Hinf$ objective. 

Unlike the pure $\Htwo$ and pure $\Hinf$ controllers, it was proved by \cite{megretski1994order} that the optimal mixed $\Htwo/\Hinf$ controller is non-rational whenever the $\Hinf$ constraint is active. Since non-rational functions do not admit finite-dimensional state-space realizations, future works focused either on finite-dimensional approximations or more tractable auxiliary objectives.

Among these, \cite{halder1998design} addressed this problem by converting it into a series of convex sub-problems. \cite{whorton1996homotopy} proposed a homotopy algorithm for fixed-order controller synthesis whereas \cite{schomig1995mixed} considered a gradient-based method. \cite{hindi1998multiobjective} proposed an approximate finite-dimensional parametrization by choosing a finite basis for the Youla parameter to compute the suboptimal solutions. The performance characteristics of the mixed $\Htwo/\Hinf$ objective have also been studied in \cite{hassibi1998upper} where the authors provide a relation between the $\Htwo$-norms of the pure $\Htwo$-optimal controller and the optimal mixed $\Htwo/\Hinf$ controller.

\subsection{Contributions}
In this work, we study the infinite-horizon mixed $\Htwo$-optimal control with $\Hinf$-norm constraints of finite-order, discrete-time, linear time-invariant (LTI) systems. While several past works \cite{hassibi1998upper, doyle1989optimal} assumed separation of the disturbances into stochastic and deterministic components, we make no such separation assumption. In particular, we consider the full-information setting where the controller can access current and past disturbances. Our major contributions are as follows:
\begin{enumerate}
    \item[i.] \textit{The Exact Stabilizing Optimal Controller:}  Moving away from earlier approaches that utilized approximation strategies and auxiliary problems to obtain a finite-dimensional optimization problem, we find the \emph{exact stabilizing optimal controller} in the frequency-domain for the infinite-horizon mixed $\Htwo/\Hinf$ problem. 
    \item[ii.] \textit{A Finite-Dimensional Characterization:} While we confirm that the optimal controller is irrational as shown in \cite{megretski1994order}, we show that a finite-dimensional parameter completely characterizes it. 
    \item[iii.] \textit{An Efficient Numerical Method:} Exploiting the finite-dimensional characterization, we propose an iterative fixed point method to compute the optimal mixed $\Htwo/\Hinf$ controller in the frequency domain to arbitrary fidelity.
    \item[iv.] \textit{Fixed-order Rational Approximations:} Given a finite order, we find the best rational approximation (in $H_{\infty}$ norm) to the optimal mixed $\Htwo/\Hinf$ controller and obtain a state space structure for the approximate controller.  We provide numerical simulations to analyze the performance of the mixed $\Htwo/\Hinf$ controller and its rational approximations for different orders.
\end{enumerate}

\section{Preliminaries}
 \textit{Notations:}
Calligraphic symbols ($\K$, $\M$, $\L$, etc.) represent operators, with $\I$ symbolizing the identity operator. Letters in boldface ($\x$, $\u$, $\w$, etc.) refer to infinite sequences. The notation $\M^\ast$ indicates the adjoint of the operator $\M$, while $\psdgeq$ signifies the Loewner order. $\tr$ denotes the normalized trace over block Laurent operators, and $\Tr$ denotes the usual trace over finite dimensional matrices. For $p\in[1,\infty]$, we denote by $\Lp_p^{n\times m}$ the Banach space of 
$n\times m$-block Laurent operators whose transfer matrix has finite $L_p$ norm on the unit circle. Similarly, $\Hp_p^{n\times m}$ is the Banach space of causal $n\times m$-block Laurent operators whose transfer matrix has finite $L_p$ norm on the unit circle. Norms $\norm[\op]{\cdot}$ and $\norm[\HS]{\cdot}$ represent the $\Hinf$/$L_\infty$ and $\Htwo$/$L_\infty$ norms, respectively.  $\norm{\cdot}$ is reserved for Euclidean norm for vectors and $\ell_2$ norm of sequences. The expressions $\{\M\}_{\!+}$ and $\{\M\}_{\!-}$ delineate the causal and strictly anti-causal portions of an operator $\M$. The absolute value of an operator is defined as $\abs{\M} \defeq \sqrt{\M^\ast \M}$. $\mathscr{T}_{+}^{m}$ denotes the set of positive symmetric polynomials of degree $m$. $\Sym_{+}^{m}$ denotes the set of symmetric positive semidefinite matrices of size $m\times m$.  $\overline{\sigma}(M)$ denotes the maximum singular value of a matrix $M$. 

\subsection{Linear-Quadratic Control}
We consider a discrete-time linear time-invariant (LTI) dynamical system described by its state-space representation as:
\begin{equation}\label{eq:state space}
x_{t+1} = A x_{t} + B_u u_{t} + B_w w_{t},
\end{equation}
Here, the vectors $x_{t} \in \R^{d_x}$, $u_{t} \in \R^{d_u}$, and $w_{t} \in \R^{d_w}$ denote the state, control input, and exogenous disturbance at time $t\in \Z$, respectively. At a given time instance $t\in\Z$, the controller suffers a per stage cost $c_t \defeq x_t^\tp Q x_t + u_t^\tp R u_t$, where $Q,R\psdg 0$. We assume that $(A,B_u)$ and $(A,B_w)$ are stabilizable. For notational convenience, we take $Q = I_{d_x}$ and $R=I_{d_u}$ without loss of generality by change of variables $x_t \mapsto Q_t^{-1/2}x_t$ and $u_t\mapsto R_t^{-1/2}u_t$ so that $c_t = \norm{x_t}^2 + \norm{u_t}^2$.

\paragraph{Input-Output Formalism}
Throughout this paper, we employ an operator-theoretic framework to represent the state-space dynamics in \cref{eq:state space} as input-output maps. We use $\x \defeq \{x_t\}_{t\in \Z}$, $\u\defeq \{u_t\}_{t\in \Z}$, and $\w\defeq \{w_t\}_{t\in \Z}$ to represent the state, control input, and disturbance sequences, respectively. The dynamical relations among these variables dictated by the state-space structure in \eqref{eq:state space} can be captured equivalently and succinctly using input/output transfer operators as:
\begin{equation}\label{eq:operator}
    \x = \F \u + \G \w,
\end{equation}
where $\F$ and $\G$ symbolize strictly causal $d_x\times d_u$ and $d_x\times d_w$-block Laurent operators mapping control inputs $\u$ and the disturbances $\w$ to the states $\x$. The corresponding transfer matrices are given by 
\begin{equation}\label{eq:transfer functions}
    F(z) \= C(zI\-A)^\inv B_u,\quad G(z) \= C(zI\-A)^\inv B_w.
\end{equation}

\paragraph{Controller} 
We consider the full-information setting where the control input $u_t$ at any time $t\in\Z$ has causal access to past disturbances, $\{w_{\tau}\}_{\tau \leq t}$. In particular, we restrict our attention to causal LTI controllers such that $u_t = \suml_{s\leq t} \widehat{K}_{t-s} w_s$ where $\{\widehat{K}_t\}_{t\geq 0}$ is the Markov parameters of the controller. We succinctly capture this relationship via a linear mapping $$\K:\w\mapsto \u \defeq \K \w$$ where $\K$ is a causal $d_u \times d_w$-block Laurent operator. Furthermore, we seek controllers that map $\ell_2$ disturbances to $\ell_2$ control inputs, which makes them bounded. Therefore, $\K$ must be a member of $\Hp_\infty^{d_u\times d_w}$.


When the underlying system in \eqref{eq:operator} is in a feedback loop with a fixed controller $\K$, the states and the control inputs are determined completely by the disturbances through the closed-loop transfer operator given by:
\begin{equation} \label{eq:T_K}
    \T_\K: \w \mapsto \begin{bmatrix} \x \\ \u \end{bmatrix} \defeq \begin{bmatrix} \F \K + \G \\ \K \end{bmatrix}\w.
\end{equation}
We use $K(z)$ and $T_K(z)$ to represent the z-domain transfer matrices corresponding to the operators $\K$ and $\T_\K$, respectively.

In the remainder of this paper, our focus will primarily be on the quadratic expression $\T_\K^\ast \T_\K \psdgeq 0$, which can be expressed in terms of $\K$ after a completion-of-squares as
\begin{equation}\label{eq:completion of squares}
    \T_{\K}^\ast \T_{\K} = (\K \- \K_\nc)^\ast (\I \+ \F^\ast \F)(\K \- \K_\nc) + \T_{\K_\nc}^\ast \T_{\K_\nc}.
\end{equation}
Here, $\K_\nc \defeq -(\I \+ \F^\ast \F)^\inv \F^\ast \G $ denotes the unique non-causal controller such that $\T_{\K_\nc}^\ast \T_{\K_\nc} \psdleq \T_{\K}^\ast \T_{\K}$ for all $\K$, \cite{hassibi_indefinite-quadratic_1999}.

\subsection{Mixed $\Htwo/\Hinf$ Control}
This paper explores the problem of mixed $\Htwo/\Hinf$ control, aiming to design a causal controller minimizing the $\Htwo$ norm of the closed-loop transfer operator, $\T_\K$, while ensuring that its $\Hinf$ norm is bounded above by a fixed constant $\gamma>0$. This can be stated formally as a constrained optimization problem as follows: 
\begin{problem}[Mixed $\Htwo/\Hinf$ Control]\label{prob:MHH}
    Given an achievable $\Hinf$ norm bound $\gamma>0$, find a causal and bounded controller {$\K\in \Hp_\infty^{d_u\times d_w}$} that achieves the minimum closed-loop $\Htwo$ norm among all $\gamma$-suboptimal $\Hinf$ controllers, \ie,
    \begin{equation} \label{eq:MHH}
        \inf_{\K\in \Hp_\infty^{d_u\times d_w}} \norm[\HS]{\T_{\K}}^2 \quad \subjto \quad \norm[\op]{\T_{\K}} \leq \gamma.
    \end{equation}
\end{problem}
We will drop the superscript over $\Lp$ and $\Hp$ spaces whenever the block dimensions are clear. Here, the $\Htwo$ norm of $\T_\K$ is defined as 
\begin{equation}
    \norm[\HS]{\T_\K} \defeq  \sqrt{\tr(\T_\K^\ast \T_\K)} = \sqrt{\int_{0}^{2\pi} \frac{d\omega}{2\pi} \Tr(T_K(\ejw)^\ast T_K(\ejw)) }
\end{equation}
The $\Htwo$ criterion can be derived from the infinite-horizon average expected cost when the disturbances constitute a white random process with $\E[w_s w_{t}^\tp ] = I_{d_w} \delta_{s-t}$,  namely,
\begin{equation}
        \limsup_{T\to \infty}  \frac{1}{2T\+1} \E\br{ \suml_{t=-T}^{T}  \norm{x_t}^2 + \norm{u_t}^2 } =\norm[\HS]{\T_\K}^2.
\end{equation}
We denote by $\Ktwo \defeq \argmin \{ \norm[\HS]{\T_{\K}} \mid \K\in \Hp_\infty^{d_u\times d_w}\}$ the optimal $\Htwo$ controller.

Similarly, the $\Hinf$ norm of $\T_\K$ corresponds to its operator norm as a mapping from $\ell_2(\Z)$ to $\ell_2(\Z)$, \ie, 
\begin{equation}
    \norm[\op]{\T_\K} \defeq \sup_{\w \in \ell_2/ \{0\}} \frac{\norm{\T_\K \w}}{\norm{\w}} = \max_{\omega \in [0,2\pi)} \overline{\sigma} (T_K(\ejw)) 
\end{equation}
Since $\T_\K$ consists of $2\times 1$ block of causal Laurent operators, the corresponding transfer matrix $T_K(z)$ is analytic outside the unit circle whenever it is bounded on the unit circle. The $\Hinf$ criterion can be derived from the worst-case infinite-horizon cost incurred by $\K$ among all bounded energy (or bounded power) disturbances, namely,
\begin{equation*}
    \sup_{\norm{\w}\leq 1} \sum_{t=-\infty}^{\infty} \norm{x_t}^2 \+ \norm{u_t}^2 = \sup_{\norm{\w}\leq 1} \norm{\x}^2 \+ \norm{\u}^2 = \norm[\op]{\T_\K}^2. 
\end{equation*}
We denote by $\Kinf \defeq \argmin\{ \norm[\op]{\T_{\K}} \mid \K\in \Hp_\infty^{d_u\times d_w}\}$ the optimal $\Hinf$ controller.

If it exists, the optimal mixed $\Htwo/\Hinf$ controller achieves the minimum expected cost for stochastic disturbances while guaranteeing a certain degree of robustness against adversarial bounded energy (or bounded power) disturbances.

Let $\ght \defeq \norm[\infty]{\T_{\Ktwo}}$ and $\ghi \defeq \norm[\infty]{\T_{\Kinf}}$ be the $\Hinf$ norms of closed-loop transfer operators under the optimal $\Htwo$ and $\Hinf$ controllers, respectively. Note that for $\gamma \geq \ght$, the $\Hinf$-norm constraint in Problem \ref{prob:MHH} is redundant, and the optimal solution coincides with $\Ktwo$. Moreover, for $\gamma < \ghi$, Problem \ref{prob:MHH} is not feasible since no causal controller can achieve $\Hinf$ norm less than that of the $\Hinf$ controller, $\Kinf$. For $\gamma = \ghi$, the solution coincides with $\Kinf$. Thus, we are interested in non-trivial solutions for $\gamma \in (\ghi, \ght)$, which interpolate between the optimal $\Htwo$ and $\Hinf$ controllers.

\section{Main Results}
In this section, we present the main theoretical results of our paper. In Theorem \ref{th:strong_duality}, we formulate the Lagrange dual of Problem \ref{prob:MHH} and establish strong duality. In Theorem \ref{thm:main}, we state the optimal controller and argue that it is irrational.



First, we form the Lagrange dual of Problem \ref{prob:MHH} in the following theorem.
\begin{theorem}[Strong Duality]
\label{th:strong_duality}
    Let $\gamma \in (\ghi, \ght)$ be an admissible $\Hinf$ norm bound. The infinite-horizon mixed $\Htwo/\Hinf$ control problem (\cref{prob:MHH}) is equivalent to the following max-min problem:
    \begin{equation}
        \max_{\substack{\Uplambda \in  \Lp_1,\\ \Uplambda \psdgeq 0}} \, \min_{\K\in \Hp_\infty} \,  \tr( \T_{\K}^\ast \T_{\K} (\I + \Uplambda)) - \gamma^2 \tr(\Uplambda),
        \label{eq:dual_MHH}
    \end{equation}
    where the dual variable  $\Uplambda \in \Lp_1^{d_w\times d_w}$ is a positive definite, self-adjoint, $d_w\times d_w$-block Laurent operator. 
\end{theorem}
\begin{IEEEproof}
    The proof of this theorem relies on the Lagrange duality theory for infinite-dimensional optimization on Banach spaces.
    Consider Problem \ref{prob:MHH}. The $\Hinf$ norm constraint $\norm[\infty]{\T_\K}\leq \gamma$ is equivalent to the $L_\infty$ norm constraint $\norm[\infty]{\T_\K^\ast \T_\K}\leq \gamma^2$, which in turn is equivalent to the  following operator inequality constraint:
    \begin{equation}
        \T_{\K}^\ast \T_{\K} - \gamma^2 \I \psdleq 0 .
    \end{equation}
    Since $\T_\K^\ast \T_\K\in \Lp_\infty^{d_w\times d_w}$ is self-adjoint, this operator inequality constraint is, by definition, equivalent to
    \begin{equation}
        \tr(\Uplambda (\T_{\K}^\ast \T_{\K} - \gamma^2 \I)) \leq 0, \quad \forall \Uplambda \psdgeq 0,
    \end{equation}
    where $\Uplambda \in \Lp_1^{d_w\times d_w}$ is a positive-definite, self-adjoint, block Laurent operator with bounded absolute trace, \ie, $\tr(\abs{\Uplambda}) < +\infty$.
     This constraint can be reincorporated into the primal problem in \eqref{eq:MHH} via Lagrangian, which yields an equivalent min-max problem:
    \begin{equation}
    \label{eq:primal_MHH_LMI}
        \inf_{\K \in \causal } \sup_{\Uplambda\psdgeq 0 } \tr( \T_{\K}^\ast \T_{\K}) + \tr(\Uplambda (\T_{\K}^\ast \T_{\K} - \gamma^2 \I)).
    \end{equation}
    Notice that the Lagrangian objective function is strictly convex (indeed quadratic) in $\K$, and affine in the dual variable $\Uplambda\psdgeq 0$. Furthermore, \emph{Slater's conditions} are satisfied since $\Kinf \in \Hp_{\infty}$ and $\norm[\op]{\T_{\Kinf}}\! \eqdef\! \ghi\< \gamma$. Thus, by \cite[Thm. 2.9.2]{zalinescu_convex_2002}, strong duality holds, and the supremum is attained. Moreover, the inner infimum is attained for a fixed $\Lambda\psdgeq 0$ due to strict convexity wrt $\K$. Rearranging terms, we get \eqref{eq:dual_MHH}.
\end{IEEEproof}
As opposed to the $\Hinf$-norm constrained primal problem in \eqref{eq:MHH}, the max-min problem in \eqref{eq:dual_MHH} is more manageable as it is not norm constrained and strictly convex in $\K$. This comes at the expense of an additional maximization step. 

Notice that the inner minimization over causal $\K$ in \eqref{eq:dual_MHH} is nothing but $(\I+\Uplambda)$-weighted squared $\Htwo$-norm objective. Therefore, the inner minimization can be considered a stochastic optimal control problem under weakly stationary disturbances with autocovariance operator $\I+\Lambda$. This can be carried out tractably using the Wiener-Hopf method \cite{youla_modern_1976} and canonical spectral factorization of $\I+\Uplambda$ as stated in \cref{lem:wiener hopf}.
\begin{lemma}[{Wiener-Hopf Method}]\label{lem:wiener hopf}
     Let $\Uplambda\in \Lp_1$ be a positive-definite and self-adjoint block Laurent operator. Consider the problem of finding an optimal causal controller, $\K$, minimizing  $(\I+\Uplambda)$-weighted $\Htwo$ norm, \ie,
    \begin{equation}\label{eq:wiener hopf}
        \min_{\K\in \Hp_\infty} \tr(\T_\K^\ast \T_\K (\I+\Uplambda)) 
    \end{equation}
    The unique solution to this problem is given by, 
    \begin{equation}
        \K_{\Uplambda} = \Ktwo + \Delta^{-1}\cl{\{\Delta \K_\nc\}_{\!-} \L}_{\!+} \L^{-1},
        \label{eq::nehari_controller}
    \end{equation}
    where $\Delta$ and $\L$ are the unique causal and causally invertible canonical spectral factors such that $\Delta^\ast \Delta = \I + \F^\ast \F$ and $\L\L^\ast = \I+\Uplambda$.
\end{lemma}
\begin{IEEEproof}
Inserting the equation \eqref{eq:completion of squares} in place of $\T_\K^\ast \T_\K$ and replacing  $\I + \F^\ast \F$ and $\I+\Uplambda$ by their corresponding canonical spectral factors, we rewrite \eqref{eq::nehari_controller} as
\begin{equation}\label{eq:nehari 1}
    \min_{\K\in \Hp_\infty } \norm[\HS]{\Delta \K \L - \Delta \K_\nc \L}^2 + \tr(\T_{\K_\nc}^\ast \T_{\K_\nc}(\I + \Uplambda)).
\end{equation}
Notice that the term $\tr(\T_{\K_\nc}^\ast \T_{\K_\nc}(\I + \Uplambda))$ does not depend on $\K$ and, therefore, can be omitted. Inside the $\Htwo$ norm, $\Delta \K \L$ is causal for any $\K$ whereas $\Delta \K_\nc \L$ is mixed (non-causal). Therefore, $\Delta \K \L$ can only minimize the objective by matching itself to the causal part of $\Delta \K_\nc \L$, \cite{youla_modern_1976, hassibi_indefinite-quadratic_1999}. The optimal solution satisfies $\Delta \K_{\Uplambda} \L = \{\Delta \K_\nc \L\}_{\!+}$ which yields $ \K_{\Uplambda} = \Delta^{-1}\cl{\Delta \K_\nc \L}_{\!+} \L^{-1}$. By splitting $\Delta\K_\nc = \{\Delta \K_\nc\}_{\!+} + \{\Delta \K_\nc\}_{\!-} $ and noting that $\Delta^{-1} \{\Delta \K_\nc\}_{\!+}$ is the $\Htwo$-controller \cite{hassibi_indefinite-quadratic_1999}, we get the desired result.
\end{IEEEproof}

\begin{remark}
    Notice that $\Uplambda=0$ yields the $\Htwo$-optimal solution since $\L\L^\ast = \I$ while any other choice of $\Uplambda\psdgeq 0$ yields a controller which is a combination of the $\Htwo$-controller, $\Ktwo$, and a compensation term $\Delta^{-1}\cl{\{\Delta \K_\nc\}_{\!-} \L}_{\!+} \L^{-1}$. This term accounts for the correlations in the disturbance as dictated by the strictly positive autocovariance operator $\I+\Uplambda$.
\end{remark}


Given any $\Uplambda\psdgeq 0$, the inner minimization solution of the Lagrange dual max-min problem in \eqref{eq:dual_MHH} can be derived in closed form using \eqref{eq::nehari_controller}. This allows us to simplify the primal problem over causal operators $\K$ in \eqref{eq:MHH}, into a dual problem over positive operators $\Uplambda\psdgeq 0$. In Theorem \ref{thm:main}, we state the optimality conditions for the dual variable $\Uplambda\psdgeq 0$.

\begin{theorem}[Saddle Point]\label{thm:main}
    Let $\gamma \inn (\ght, \ghi)$ be a fixed admissible $\Hinf$ norm bound and $(\K_{\gamma}, \Uplambda_{\gamma})$ be a saddle point of the max-min problem in \eqref{eq:dual_MHH}. Moreover, let $\L_\gamma$ and $\Delta$ be the unique causal and causally invertible spectral factors of $\I\+\Uplambda_{\gamma} \= \L_{\gamma} \L_{\gamma}^\ast$ and $\I\+\F^\ast \F \= \Delta^\ast \Delta$, respectively. Then, $(\K_{\gamma}, \Uplambda_{\gamma})$ satisfies the following necessary and sufficient conditions\footnote{$\max\{\clf{X}, \clf{Y}\} \triangleq \frac{1}{2}(\clf{X} + \clf{Y} + \abs{\clf{X} - \clf{Y}})$}:
    \begin{align}
        \K_{\gamma} &= \K_{2} + \Delta^{-1}\cl{\{\Delta \K_\nc\}_{\!-} \L_{\gamma}}_{\!+} \L_{\gamma}^{-1}, \quad \text{and} \label{eq:controller_opt} \\
        \L_{\gamma}^{\ast} \L_{\gamma} &= \max\left\{\gamma^{-2}\pr{\Ss^\ast\Ss +  \L_{\gamma}^\ast \T_{\K_{\nc}}^{\ast}\T_{\K_{\nc}}\L_{\gamma}},\; \I\right\},
                \label{eq:M_opt_MIMO}
    \end{align}
    where $\Ss_{\gamma} \defeq \{\Delta\K_{\nc}\L_{\gamma}\}_{\!-}$. 
    \label{th:saddle}
\end{theorem}
\begin{corollary}
    When the disturbances are scalar ($d_w \= 1$), we have,
    \begin{align}
        &\L_{\gamma}^\ast \L_{\gamma} = \max \left\{ \frac{\{\Delta\K_{\nc}\L_{\gamma}\}_{\!-}^\ast \{\Delta\K_{\nc}\L_{\gamma}\}_{\!-}}{\gamma^2 - \T_{\K_{\nc}}^{\ast}\T_{\K_{\nc}}}, \;\I\right \}.
        \label{eq:M_opt}
    \end{align}
\end{corollary}
\begin{IEEEproof}
    Existence of a saddle point $(\K_{\gamma}, \Uplambda_{\gamma})$ is ensured by \cref{th:strong_duality}. We can characterize the saddle point by a Banach space analog of Karush-Kuhn-Tucker (KKT) conditions \cite[Thm. 2.9.2]{zalinescu_convex_2002}:
    \begin{itemize}
        \item[i.] Stationarity: $\K_\gamma \in \argmin\{ \tr(\T_\K^\ast \T_\K (\I+\Uplambda_\gamma))\mid {\K \in \Hp_\infty}\}$,
        \item[ii.] Primal feasibility: $ \T_{\K_\gamma}^\ast \T_{\K_\gamma} - \gamma^2 \I \psdleq 0 $,
        \item[iii.] Dual feasibility: $\Uplambda_\gamma \psdgeq 0$,
        \item[iv.] Complementary slackness: $ \tr(\Uplambda_\gamma (\T_{\K_\gamma}^\ast \T_{\K_\gamma} - \gamma^2 \I))=0$.
    \end{itemize}
    By the Wiener-Hopf method in \cref{lem:wiener hopf}, the stationarity condition (i) immediately implies Eq. \eqref{eq:controller_opt}. 

To characterize the dual variable $\Uplambda_\gamma$, we note that the complementary slackness condition (iv) can be equivalently expressed as 
\begin{equation}
    \Uplambda_{\gamma} \left( \T_{\K_{\gamma}}^\ast \T_{\K_{\gamma}} - \gamma^2 \I \right) = 0.
    \label{eq:slackness}
\end{equation}
since $ \T_{\K_\gamma}^\ast \T_{\K_\gamma} - \gamma^2 \I \psdleq 0 $ and $\Uplambda_\gamma \psdgeq 0$ by primal (ii) and dual (iii) feasibility. Furthermore, this condition can be written alternatively as $  ( \T_{\K_{\gamma}}^\ast \T_{\K_{\gamma}} - \gamma^2 \I )\Uplambda_{\gamma} = 0$. Therefore, the operators $\T_{\K_{\gamma}}^\ast \T_{\K_{\gamma}}-\gamma^2 \I$ and $\Uplambda_{\gamma}$ commute. This implies that they share the same eigenspace and can be simultaneously unitarily diagonalized, so do the operators $\T_{\K_{\gamma}}^\ast \T_{\K_{\gamma}}$ and $\Uplambda_{\gamma}$.

For the optimal controller $\K_\gamma$ in  \eqref{eq:controller_opt}, we have 
\begin{align}
    \Delta \K_\gamma-\Delta \K_\nc &= \{\Delta \K_\nc \L_\gamma\}_{+}\L_\gamma^{-1} - \Delta \K_\nc \L_\gamma \L_\gamma^{-1},\\
    &=-\{\Delta \K_\nc \L_\gamma\}_{-} \eqdef -\Ss_\gamma.
\end{align}
Therefore, the quadratic expression $ \T_{\K_{\gamma}}^\ast \T_{\K_{\gamma}}$ becomes
\begin{align}
     \T_{\K_{\gamma}}^\ast \T_{\K_{\gamma}}  &= (\Delta \K_\gamma-\Delta \K_\nc)^\ast (\Delta \K_\gamma-\Delta \K_\nc) +  \T_{\K_{\nc}}^{\ast}\T_{\K_{\nc}},\\
     &=\L_{\gamma}^{-\ast} \Ss^\ast  \Ss \L_{\gamma}^{-1}
     + \T_{\K_{\nc}}^{\ast}\T_{\K_{\nc}}.
\end{align}
Substituting the optimal value of $\T_{\K_{\gamma}}^\ast \T_{\K_{\gamma}}$ in \eqref{eq:slackness}, we get,
\begin{equation}
    \Uplambda_{\gamma} \left( \L_{\gamma}^{-\ast} \Ss^\ast \Ss \L_{\gamma}^{-1}
     + \T_{\K_{\nc}}^{\ast}\T_{\K_{\nc}} - \gamma^2 \I \right) = 0.
 \label{eq:kkt_with_l}
\end{equation}
Using the identity $\Uplambda_{\gamma} =  \L_{\gamma}\L_{\gamma}^{\ast} - \I$, we get,
\begin{equation*}
     \L_{\gamma} \Ss^\ast \Ss \L_{\gamma}^{-1}  + \left(\L_{\gamma}\L_{\gamma}^{\ast} - \I\right) \left(\T_{\K_{\nc}}^{\ast}\T_{\K_{\nc}} - \gamma^2 \I\right) - \L_{\gamma}^{-\ast} \Ss^\ast  \Ss \L_{\gamma}^{-1} = 0.
\end{equation*}
Multiplying the identity above by $\L_{\gamma}$ on the right and by $\L_{\gamma}^\ast$ on the left, we get,
\begin{equation*}
    \L_{\gamma}^\ast\L_{\gamma} \Ss^\ast  \Ss  + \L_{\gamma}^\ast\left(\L_{\gamma}\L_{\gamma}^{\ast} - \I\right)   \left(\T_{\K_{\nc}}^{\ast}\T_{\K_{\nc}} - \gamma^2 \I\right)\L_{\gamma}  -  \Ss^\ast \Ss = 0.
\end{equation*}
The identity above can be simplified further by factorizing $\L_{\gamma}^{\ast} \L_{\gamma} - \I$ as follows:
\begin{align}
     0&=(\L_{\gamma}^\ast\L_{\gamma}- \I) \Ss^\ast  \Ss  + \left(\L_{\gamma}^\ast\L_{\gamma}\L_{\gamma}^{\ast} - \L_{\gamma}^\ast \right)   \left(\T_{\K_{\nc}}^{\ast}\T_{\K_{\nc}} - \gamma^2 \I\right)\L_{\gamma},\nonumber\\
     &=(\L_{\gamma}^\ast\L_{\gamma}- \I) \Ss^\ast  \Ss  + \left(\L_{\gamma}^\ast\L_{\gamma} - \I \right) \L_{\gamma}^{\ast}   \left(\T_{\K_{\nc}}^{\ast}\T_{\K_{\nc}} - \gamma^2 \I\right)\L_{\gamma},\nonumber\\
     &=(\L_{\gamma}^\ast\L_{\gamma}- \I) \pr{\Ss^\ast  \Ss  + \L_{\gamma}^{\ast} \T_{\K_{\nc}}^{\ast}\T_{\K_{\nc}}\L_{\gamma} - \gamma^2 \L_{\gamma}^\ast \L_{\gamma}}.\label{eq:version1}
\end{align}
Starting with $ \T_{\K_\gamma}^\ast \T_{\K_\gamma} - \gamma^2 \I \psdleq 0 $ instead, a similar identity can be derived: 
\begin{equation}
     0=\pr{\Ss^\ast  \Ss  + \L_{\gamma}^{\ast} \T_{\K_{\nc}}^{\ast}\T_{\K_{\nc}}\L_{\gamma} - \gamma^2 \L_{\gamma}^\ast \L_{\gamma}}(\L_{\gamma}^\ast\L_{\gamma}- \I) .\label{eq:version2}
\end{equation}
Therefore, by \eqref{eq:version1} and \eqref{eq:version2}, we have that  $\Ss^\ast  \Ss  + \L_{\gamma}^{\ast} \T_{\K_{\nc}}^{\ast}\T_{\K_{\nc}}\L_{\gamma}$ and $\L_{\gamma}^\ast\L_{\gamma}$ commute, and thus are simultaneously diagonalized.

Defining $\clf{X}_\gamma \defeq \gamma^{-2}(\Ss^\ast  \Ss  + \L_{\gamma}^{\ast} \T_{\K_{\nc}}^{\ast}\T_{\K_{\nc}}\L_{\gamma})$ and $\clf{N}_\gamma \defeq \L_\gamma^\ast \L_\gamma$, we can rewrite \eqref{eq:version1} and \eqref{eq:version2} as
\begin{equation}\label{eq:equation to solve}
    (\NN_\gamma - \I )(\NN_\gamma -\clf{X}_\gamma )=0,
\end{equation}
with $\NN_\gamma - \I \psdgeq 0$ by dual feasibility, $\NN_\gamma -\clf{X}_\gamma\psdgeq 0$ by primal feasibility, and commuting operators $\NN_\gamma$ and $\clf{X}_\gamma$. Since these operators commute, we can solve the quadratic equation \eqref{eq:equation to solve} for $\NN_\gamma$ in terms of $\clf{X}_\gamma$. This essentially gives,
\begin{equation}
    \NN_\gamma = \half \pr{ \clf{X}_\gamma + \I + \abs{\clf{X}_\gamma - \I}} = \max\cl{\clf{X}_\gamma, \, \I }
\end{equation}
which is the only possible solution of \eqref{eq:equation to solve} that satisfies the primal and dual feasibility constraints. This expression is essentially an operator-theoretic analog of taking the maximum between two elements. Thus, we get \eqref{eq:M_opt_MIMO}.

When the disturbances are scalar, \ie, $d_w = 1$, all $d_w\times d_w$-block Laurent operators commute with each other. Thus, the expression \eqref{eq:version2} can be simplified further as
\begin{equation}
     0=\pr{\Ss^\ast  \Ss  + \L_{\gamma}^{\ast}\L_{\gamma}(\T_{\K_{\nc}}^{\ast}\T_{\K_{\nc}} - \gamma^2)  }(\L_{\gamma}^\ast\L_{\gamma}- \I).
\end{equation}
Upon solving it for $\L_\gamma^\ast \L_\gamma$ in the same fashion, we immediately get \eqref{eq:M_opt}

\end{IEEEproof}

\section{A Fixed-Point Iteration Algorithm}
Henceforth in this paper, we focus our analysis on the case when $d_w = 1$, \ie, scalar disturbances.
In subsection \ref{subsec::subK}, we first demonstrate that the Karush-Kuhn-Tucker (KKT) conditions outlined in equation~\eqref{eq:M_opt} can be exclusively defined by a finite-dimensional parameter, $\overline{B}_{\gamma}$, within the frequency domain. However, the optimal controller lacks a rational nature, preventing it from being represented through a finite-dimensional state-space model. Following this, in subsection \ref{subsec::alg}, we introduce a fixed-point iteration method for any given $\gamma \in (\ghi, \ght)$. This method is designed to determine $\overline{B}_{\gamma}$, enabling the calculation of the optimal controller, $K_{\gamma}(\e^{j\omega})$ in the frequency domain.

\subsection{Finite-Dimensional Parameterization of the Optimal Controller}\label{subsec::subK}
In the following Theorem~\ref{lemma:finiteBl}, we establish that the strictly anticausal transfer function $S_{\gamma,-}(\e^{j\omega})$ defined below can be expressed using a finite-dimensional state-space model.
This theorem indicates that the right-hand side of \eqref{eq:N_freq} for $N_\gamma(\e^{j\omega})$, involving the square root of the rational term $S_{\gamma,-}(\e^{j\omega})^\ast S_{\gamma,-}(\e^{j\omega})$, is no longer rational due to the square root operation. This observation leads us to Corollary~\ref{thm:irrational}.

\begin{theorem}\label{lemma:finiteBl}
    Define $S_{\gamma,\-}(\e^{j\omega}) \defeq \{\Delta K_\circ L_\gamma\}_{\-}(\e^{j\omega})$. Then,
    \begin{align}
    &S_{\gamma,-}(\e^{j\omega}) = \overline{C}(\e^{-j\omega}I - \overline{A})^{-1} \overline{B}_{\gamma}, \label{eq:sg_freq}\\
    &\quad \text{with} \quad \overline{B}_{\gamma} = \frac{1}{2\pi} \int_{0}^{2\pi} (I-\e^{j\omega} \overline{A})^{-1} \overline{D} L_\gamma(\e^{j\omega}) d\omega.
    \end{align}
    The optimal controller in the frequency domain is given by,
    \begin{align}
        K_\gamma(\e^{j\omega}) &= K_{\nc}(\ejw)\- \Delta^\inv (\e^{j\omega}) S_{\gamma,\-}(\e^{j\omega}) L_\gamma^\inv (\e^{j\omega}), \label{eq:K_freq}\\
    N_\gamma(\e^{j\omega}) &=  \max\left\{\frac{\{\Delta\K_{\nc}\L_{\gamma}\}^{\ast}_{\!-}(\ejw)\{\Delta\K_{\nc}\L_{\gamma}\}_{\!-}(\ejw)}{\gamma^2 - T_{\K_{\nc}}(\ejw)^{\ast}T_{\K_{\nc}}(\ejw)} , 1\right\}\label{eq:N_freq},
    \end{align}
    where,
    \begin{align}
        &K_{\textrm{lqr}} \defeq {(R + B_u^\top PB_u)}^{-1}B_u^\top PA \\
        & P = Q + A^\top PA - A^\top PB_u{(R + B_u^\top PB_u)}^{-1}B_u^\top PA \\
        &\overline{A} \defeq \left(A - B_u K_{\textrm{lqr}}\right)^\top \\
        &\overline{C} \defeq -{(R + B_u^\top P B_u)}^{-\top/2}B_u^\top \\
        &\overline{D} \defeq \left(A - B_u K_{\textrm{lqr}}\right)^\top PB_w \\
        & N_\gamma(\e^{j\omega}) \defeq {L_\gamma(\e^{j\omega})}^\ast L_\gamma(\e^{j\omega}).
    \end{align}
\end{theorem}
\begin{IEEEproof}
Using the identity $\{\clf{X}\}_{\+}\=\clf{X}\-\{\clf{X}\}_{\-}$ , we restate the KKT equations~\eqref{eq:M_opt} in the frequency domain as in \eqref{eq:K_freq}, \eqref{eq:N_freq}.

We introduce the Linear Quadratic Regulator (LQR) controller $K_{\text{lqr}}$, the corresponding closed-loop matrix $A_K$, and the unique solution to the LQR Riccati equation that stabilizes the system, $P \succ 0$ to write $S_{\gamma,-}(\e^{j\omega})$ in \eqref{eq:sg_freq}. The Riccati equation emerges from the spectral factorization of $\Delta^\ast \Delta = \I + \F^\ast \F$. See \cite[Lemma 11]{kargin2023wasserstein} for example.
\end{IEEEproof}

\begin{corollary}\label{thm:irrational}
For any given $\gamma \in (\ghi, \ght)$, both $N_\gamma(\e^{j\omega})$ and the optimal mixed $\Htwo/\Hinf$ controller $K_\gamma(\e^{j\omega})$, are characterized as irrational. Consequently, $K_\gamma(\e^{j\omega})$ cannot be realized by a finite-dimensional state-space model.
\end{corollary}
Although $K_\gamma(\e^{j\omega})$ cannot be modeled in a finite-dimensional state-space, Lemma~\ref{lemma:finiteBl} introduces a finite-dimensional parameterization for $N_\gamma(\e^{j\omega})$ via $\overline{B}_{\gamma}$. Theorem~\ref{thm:fixed-point} further verifies that $\overline{B}_{\gamma}$ directly determines $N_\gamma(\e^{j\omega})$, and by extension, the suboptimal controller $K_\gamma(\e^{j\omega})$.
\begin{theorem}[Fixed-Point Solution]\label{thm:fixed-point}
Assuming $d_w=1$, $\gamma \in (\ghi, \ght)$, we consider a sequence of mappings:
\begin{align}
&F_1: \overline{B} \mapsto S_{-}(\e^{j\omega}) = \overline{C}(\e^{-j\omega}I - \overline{A})^{-1} \overline{B} \\
&F_{2,\gamma}: S_{-}(\e^{j\omega})  \mapsto N_\gamma(\e^{j\omega}) \nonumber \\
& \qquad=  \max\left\{\frac{S_{-}(\e^{j\omega})^\ast S_{-}(\e^{j\omega})}{\gamma^2 - T_{\K_{\nc}}(\ejw)^{\ast}T_{\K_{\nc}}(\ejw)} , 1\right\}. \label{eq:optimal_n} \\
&F_3: N(\e^{j\omega})  \mapsto L(\e^{j\omega}) \\
&F_4: L(\e^{j\omega})  \mapsto \overline{B} = \frac{1}{2\pi} \int_{0}^{2\pi} (I-\e^{j\omega} \overline{A})^{-1} \overline{D} L(\e^{j\omega}) d\omega.
\end{align}
where $F_3$ produces a unique spectral factor of $N(\e^{j\omega})>0$. The composite function $F_4 \circ F_3 \circ F_{2,\gamma} \circ F_1: \overline{B}  \mapsto \overline{B}$ has a unique fixed point $\overline{B}_{\gamma}$, with $N{\gamma}(\e^{j\omega}) \equiv F_{2,\gamma} \circ F_1 (\overline{B}_{\gamma})$ fulfilling the KKT conditions~\eqref{eq:M_opt}.
\end{theorem}
\begin{IEEEproof}
    Consider the optimality condition \eqref{eq:N_freq}. Note that for $\gamma > \ghi$, $N_\gamma(\e^{j\omega})$ in \eqref{eq:N_freq} is well-defined. Thus, the map $F_4\circ F_3\circ F_{2_{\gamma}}\circ F_1:\bar{B} \mapsto \bar{B}$ described earlier admits a fixed point $\overline{B}_\gamma$ for a fixed $\gamma$. Since \eqref{eq:primal_MHH_LMI} is concave in $\M$ (or $\NN = \M + \I$), the optimal solution $\NN^{*}$ is unique. Given that $M_\gamma(\e^{j\omega})=L_\gamma(\e^{j\omega})L_\gamma^\ast(\e^{j\omega})$, where $L_\gamma(\e^{j\omega})$ represents a spectral factor of $M_\gamma(\e^{j\omega})$ that is both causal and causally invertible, it follows that $L_\gamma(\e^{j\omega})$ is uniquely determined apart from a unitary transformation. By establishing a specific choice for the unitary transformation during the spectral factorization process, for example, opting for positive-definite factors at infinity as outlined by \cite{ephremidze_algorithmization_2018}, we ensure the uniqueness of $L_\gamma(\e^{j\omega})$. Consequently, with $\bar{A}$ and $\bar{D}$ being constants, the expression for $\bar{B}_\gamma=\frac{1}{2\pi} \int_0^{2\pi} (I-\e^{j\omega} \overline{A})^{-1} \overline{D} L_\gamma(\e^{j\omega}) d\omega$ is also uniquely determined.
\end{IEEEproof}

\subsection{Algorithm}\label{subsec::alg}
Now that we know a fixed point solution exists to our problem, we can use Theorem \ref{thm:fixed-point} to obtain the following iterative fixed-point algorithm. Once we have an optimal $N_{\gamma}(e^{j \omega})$, we can find the optimal mixed $\Htwo/\Hinf$ controller $K_{\gamma}(\e^{j \omega})$ using \eqref{eq:K_freq}.
\begin{algorithm}
   \caption{\FixAlg:  Dual Fixed-Point Iteration via Spectral Factorization  }
   \label{alg:fixed_point}
\begin{algorithmic}
   \STATE {\bfseries Input:} $\gamma\>{\gamma}_{\,H_\infty}$, initialise $\bar{B}$, system $(\overline{A}, \overline{D}, \overline{C})$ 
   \REPEAT 
   \STATE $N_\gamma^{(n)}(\ejw) \leftarrow   F_{2, \gamma}\!\circ\!F_1(\bar{B}^{(n)})$
   \STATE $L_\gamma^{(n+1)}(\mathrm{e}^{j\omega}) \leftarrow \texttt{SpectralFactor}( N_\gamma^{(n)}(\mathrm{e}^{j\omega}))$
   \STATE $\bar{B}^{n + 1} = F_4(L_\gamma^{(n+1)}(\mathrm{e}^{j\omega}))$
   \STATE $n\leftarrow n+1$
   \UNTIL{convergence of $N_\gamma^{(n)}(\mathrm{e}^{j\omega})$}
\end{algorithmic}
\end{algorithm}

\section{Convergence Analysis}
\label{sec::proof_conv}
In this section, we provide a proof of convergence of \FixAlg for the particular case when $d_x = d_w = 1$. We first show the fixed point, that Algorithm \ref{alg:fixed_point} converges to, is unique when $\gamma > \ghi$. Then, we show that the iterates produce a sequence of monotonically increasing spectra $\NN$. Due to this, the algorithm must converge to the unique fixed point. Formally, we state this in the following lemma.

\begin{lemma}[Monotonicity]\label{lem:monotone}
    For $d_x=d_w=1$, the composite mapping $F_{2,\gamma} \circ F_1 \circ F_4 \circ F_3:N(\ejw)\mapsto N_{\gamma}(\ejw)$ is monotonic, \ie, for any two positive power spectral densities such that $0 < N_1(e^{j \omega}) \leq N_2(e^{j \omega})$  for all $\omega \in [0,2\pi)$, the mapping $F_{2,\gamma} \circ F_1 \circ F_4 \circ F_3$ preserves their order.
\end{lemma}
\begin{IEEEproof}
Let $\gamma > \ghi$. Consider now two spectra $0\psdl \NN_1 \preccurlyeq \NN_2$ that are represented in the frequency domain as $N_1(e^{j \omega}) \leq N_2(e^{j \omega})$  for all $\omega \in [0,2\pi)$. Now the spectra $\NN_1, \NN_2$ are passed through one iteration of Algorithm \ref{alg:fixed_point}, \ie, $F_3\circ F_4\circ F_1 \circ F_{2_{\gamma}}\circ :\NN \mapsto \bar{\NN}$ to get $\bar{\NN_1}, \bar{\NN_2}$ respectively. We want to show that $\bar{n}_1(e^{j \omega}) \leq \bar{n}_2(e^{j \omega}) \forall \omega$ \ie, Algorithm \ref{alg:fixed_point} preserves the order of $\NN_1 \preccurlyeq \NN_2$. We have that
\begin{align}
    \bar{N_1}(\e^{j\omega}) &= \max \left\{ 1,  \frac{\{\Delta\K_{\nc}\L_{1}\}^{\ast}_{\!-}(\ejw)\{\Delta\K_{\nc}\L_{1}\}_{\!-}(\ejw)}{\gamma^2 - T_{\K_{\nc}}(\ejw)^{\ast}T_{\K_{\nc}}(\ejw)}\right\} \label{eq:n1} \\
    \bar{N_2}(\e^{j\omega}) &= \max \left\{ 1,  \frac{\{\Delta\K_{\nc}\L_{2}\}^{\ast}_{\!-}(\ejw)\{\Delta\K_{\nc}\L_{2}\}_{\!-}(\ejw)}{\gamma^2 - T_{\K_{\nc}}(\ejw)^{\ast}T_{\K_{\nc}}(\ejw)}\right\}.\label{eq:n2}
\end{align}
Moreover, denoting by $\NN_1=\L_1^\ast \L_2$ and $\NN_2=\L_2^\ast \L_2$ the unique spectral factorizations, we have that
\begin{align}
&\tr \left(\{\Delta\K_{\nc}\L_{1}\}^{\ast}_{\!-}\{\Delta\K_{\nc}\L_{1}\}_{\!-}\right) \!\stackrel{(a)}{=}\! \inf_{\K \in \Hp_\infty } \tr((\T_\K^\ast \T_\K - \T_{\K_{\nc}}^\ast \T_{\K_{\nc}}) \NN_1)  \nonumber  \\
&\quad\quad\quad\stackrel{(b)}{\leq} \inf_{\K \in \Hp_\infty} \tr((\T_\K^\ast \T_\K - \T_{\K_{\nc}}^\ast \T_{\K_{\nc}})\NN_2) \label{eq:n1leqn2} \\
&\quad\quad\quad \stackrel{(c)}{=} \tr \left(\{\Delta\K_{\nc}\L_{2}\}^{\ast}_{\!-}\{\Delta\K_{\nc}\L_{2}\}_{\!-} \right), \nonumber
\end{align}
where the identities $(a)$, $(c)$ are due to Wiener-Hopf technique in \cref{lem:wiener hopf} and $(b)$ is due to $\NN_1 \preccurlyeq \NN_2$.

Notice that, for $d_x=d_w=1$, we have that  
\begin{align}
    &\{\Delta\K_{\nc}\L_{1}\}^{\ast}_{\!-}\{\Delta\K_{\nc}\L_{1}\}_{\!-} = \abs{\bar{C}}^2\abs{\bar{B}_1}^2 \abs{\e^{-j\omega}I - \overline{A}}^{-2}, \\
     &\{\Delta\K_{\nc}\L_{2}\}^{\ast}_{\!-}\{\Delta\K_{\nc}\L_{2}\}_{\!-} = \abs{\bar{C}}^2\abs{\bar{B}_2}^2 \abs{\e^{-j\omega}I - \overline{A}}^{-2}, 
\end{align} 
and therefore, their traces take the form
\begin{align}
    &\tr\pr{\{\Delta\K_{\nc}\L_{1}\}^{\ast}_{\!-}\{\Delta\K_{\nc}\L_{1}\}_{\!-}} \nonumber \\
    &\quad\quad= \abs{\bar{C}}^2\abs{\bar{B}_1}^2 \int_{0}^{2\pi}\abs{\e^{-j\omega}I - \overline{A}}^{-2} \frac{d \omega}{2\pi}, \\
    &\tr\pr{\{\Delta\K_{\nc}\L_{2}\}^{\ast}_{\!-}\{\Delta\K_{\nc}\L_{2}\}_{\!-}}\nonumber \\
&\quad\quad=\abs{\bar{C}}^2\abs{\bar{B}_2}^2\int_{0}^{2\pi} \abs{\e^{-j\omega}I - \overline{A}}^{-2}\frac{d \omega}{2\pi}.
\end{align} 
By this observation and using the inequality \eqref{eq:n1leqn2}, we have that $\abs{\bar{B}_1}^2 \leq \abs{\bar{B}_1}^2 $ and thus the monotonicity of the traces implies the monotonicity of the operators, \ie, for all $\omega \in[0,2\pi)$, the following holds
\begin{align}
    &\{\Delta\K_{\nc}\L_{1}\}^{\ast}_{\!-}(\ejw)\{\Delta\K_{\nc}\L_{1}\}_{\!-}(\ejw)\nonumber \\
    &\quad\quad= \abs{\bar{C}}^2\abs{\bar{B}_1}^2 \abs{\e^{-j\omega}I - \overline{A}}^{-2}, \\
    &\quad\quad\leq \abs{\bar{C}}^2\abs{\bar{B}_2}^2 \abs{\e^{-j\omega}I - \overline{A}}^{-2}, \\
    &\quad\quad= \{\Delta\K_{\nc}\L_{2}\}^{\ast}_{\!-} (\ejw)\{\Delta\K_{\nc}\L_{2}\}_{\!-}(\ejw).
\end{align} 
Hence, we have that $\bar{N_1}(\e^{j\omega}) \leq \bar{N_2}(\e^{j\omega})$ for all $\omega\in [0,2\pi)$ using \eqref{eq:n1} and \eqref{eq:n2}.
\end{IEEEproof}

We can now initialise the algorithm with $\bar{B}^{(0)} = 0$ or $\bar{N}^{(0)}(\e^{j\omega}) = 1$. Since after each iteration, $\bar{N}^{(n)}(\e^{j\omega}) \geq 1$ for all $\omega\in [0,2\pi)$, \FixAlg generates a monotonically increasing sequence of $\{\bar{N}^{(n)}(\e^{j\omega})\}$ for all $\omega\in [0,2\pi)$, which converges to the unique fixed point. We state this formally in the following theorem.
\begin{theorem}
    For $d_x=d_w=1$ and $\bar{B}^{(0)}=0$, the sequence of iterates $\{\bar{N}^{(n)}(\e^{j\omega})\}$ generated by \FixAlg is monotonically increasing and converges to the optimal solution in \eqref{eq:M_opt}.
\end{theorem}
\begin{IEEEproof}
    The proof follows directly from the repeated application of monotonicity result on Lemma~\ref{lem:monotone} and \cite[Thm. 2]{yates_framework_1995}
\end{IEEEproof}

\begin{remark}
    Although we present a proof of convergence for the particular case of scalar systems $d_x=1$, empirical evidence suggests (see Section \ref{sec::numerical}) that the algorithm is exponentially convergent with a faster rate of convergence for larger $\gamma>0$. 
\end{remark}

\section{Rational Approximation}
This section describes a practical approach to devising state-space controllers that serve as approximations to our irrational $\gamma$-optimal controller \eqref{eq:K_freq}. Rather than attempting to approximate the controller directly, we \emph{approximate the power spectrum $N(\e^{j\omega})$}, aiming to reduce the $\Hinf$-norm of the approximation error through the use of positive rational functions. While the problem of approximating with rational functions typically does not lead to convex formulations, we demonstrate in Theorem \ref{thm:rational approx via feasibility} that the process of approximating positive power spectra through the use of ratios of positive fixed-order polynomials can indeed be framed as a convex feasibility problem.
Formally, the problem is stated as follows.
\begin{problem}[Rational Approximation via $\Hinf$-norm] \label{prob:hinf_rational_opt}
For a given positive spectrum $\NN$, identify the optimal rational approximation utilizing the $\Hinf$ norm, with an order not exceeding $m\in\N$. Specifically,
\begin{equation}\label{eq:hinf_rational_opt}
\inf_{\PP, \QQ \in \trig_{+}^{m}} \Norm[\op]{{\PP/\QQ} - \NN } \text{ subject to } \tr(\QQ)=1,
\end{equation}
where $\trig_{m,+}$ is the set of positive symmetric polynomials of order less than or equal to $m$ and the constraint $\Tr(\QQ)=1$ is to avoid redundancy in solutions.
\end{problem}
\begin{definition}\label{def:sublevel}
    Given an $\epsilon>0$ approximation bound, the $\epsilon$-sublevel set of Problem \ref{prob:hinf_rational_opt} is defined as
    \begin{equation*}
    \sub_\epsilon \!\defeq\! \cl{(\PP,\QQ)  \,\left|\; \Norm[\op]{{\PP/\QQ} \- \clf{N} } \!\!\leq\! \epsilon, \; \tr(\QQ)\=1 \right.} 
\end{equation*}
\end{definition}
\begin{theorem}[{Feasibility of $\sub_{\epsilon}$, \cite[Thm. 5.5]{kargin_infinite-horizon_2024}}]\label{thm:rational approx via feasibility}
    Given an accuracy level $\epsilon\>0$ and $m\in \N$ is a fixed order, the polynomials $\PP$ and $\QQ$ of order $m$ belong to the $\epsilon$-sub-level set, \ie $(\PP,\QQ) \in \sub_\epsilon$ if and only if there exists $\mathbf{P},\mathbf{Q} \in \Sym_{+}^{m +1}$ such that $\tr\pr{\mathbf{Q}} = 1$ and for all $\omega \in [0,2\pi)$, .
\begin{align}
\vspace{-3mm}
    1)\,&\Tr\pr{\mathbf{P} \bm{\Theta}(\ejw)} \- \pr{ N(\ejw)  \+ \epsilon } \Tr\pr{\mathbf{Q} \bm{\Theta}(\ejw)}\leqq 0, \\
    2)\,&\Tr\pr{\mathbf{P} \bm{\Theta}(\ejw)} \- \pr{ N(\ejw)  \- \epsilon } \Tr\pr{\mathbf{Q} \bm{\Theta}(\ejw)}\geqq 0,
    \vspace{-6mm}
\end{align}
\end{theorem}
Although the above equations hold for all frequencies, practical implementation necessitates focusing on a finite selection of frequency samples. To sidestep this limitation, the analysis can be limited to a select set of frequencies, defined as $\Omega_N = \{\omega = 2\pi k /N | k=1,\ldots,N\}$, where $N$, chosen to be much larger than $m$, provides a dense sampling of the frequency domain. While this approach inherently approximates the full spectrum of frequencies, increasing the number of sampled frequencies $N$ allows for an arbitrary improvement in the approximation's accuracy. Consequently, this method facilitates the transformation of the rational function approximation problem into a convex feasibility problem, addressable through the application of Linear Matrix Inequalities (LMIs) alongside a finite collection of affine (in)equality constraints.
Upon achieving a rational approximation of $N(\e^{j\omega})$, we then derive a state-space controller as outlined in \eqref{eq:controller_opt}. Note that a canonical factor of the rational approximation of $N(\e^{j\omega})$ can always be found due to the following lemma.

\begin{lemma}[{Canonical Factorization \cite[Lem. 1]{sayed_survey_2001}}]\label{lem:spect}
Given a Laurent polynomial of order $m$, $P(z)=\sum_{k=-m}^{m} p_k z^{-k}$, where $p_k = p_{-k} \in \R$, and $P(\e^{j\omega})> 0$, it can be shown that a canonical factor $L(z)= \ell_0 + \ell_1 z^{-1} + \ldots + \ell_m z^{-m}$ exists. This factor satisfies $P(\e^{j\omega}) = |L(\e^{j\omega})|^2$, with all roots of $L(z)$ lying inside the unit circle.
\end{lemma}


\section{Numerical Results}
\label{sec::numerical}
\begin{figure}[t]
	\centering
		\centering

\includegraphics[width=0.4\textwidth]{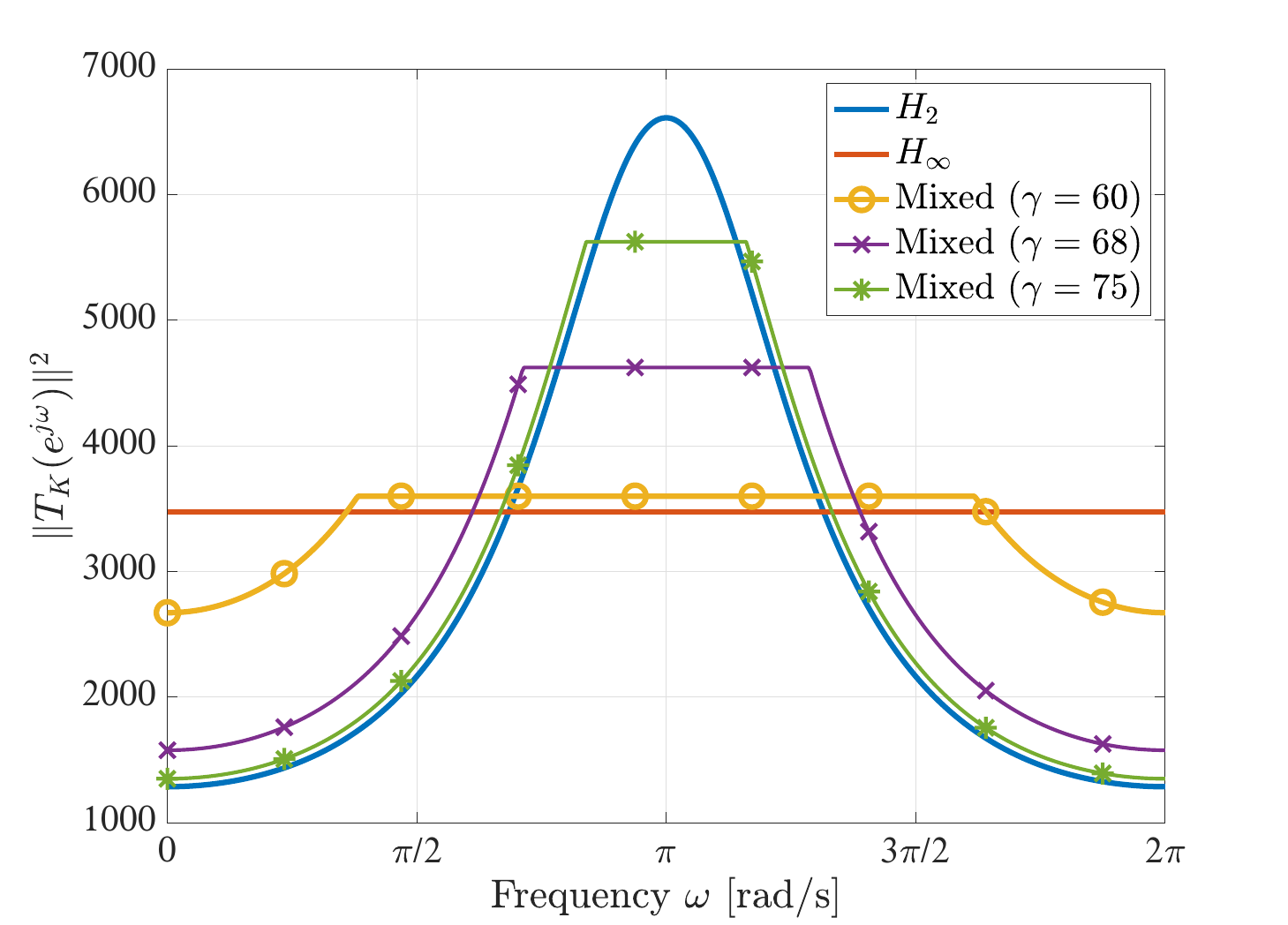}
	\caption{The spectral norm, $\overline{\sigma}(T_K^\ast(\ejw)T_K(\ejw))$ of the mixed $\Htwo/\Hinf$ controller for $\gamma \in \{60, 68, 75\}$ at different frequency values, for the system [AC17]. The cost of the mixed $\Htwo/\Hinf$ controller follows $H_2$ and clips at the threshold $\gamma$ for some frequencies.} 
 
	\label{fig:ac17}
\end{figure}

\begin{figure}[t]
	\centering
		\centering
		\includegraphics[width=0.4\textwidth]{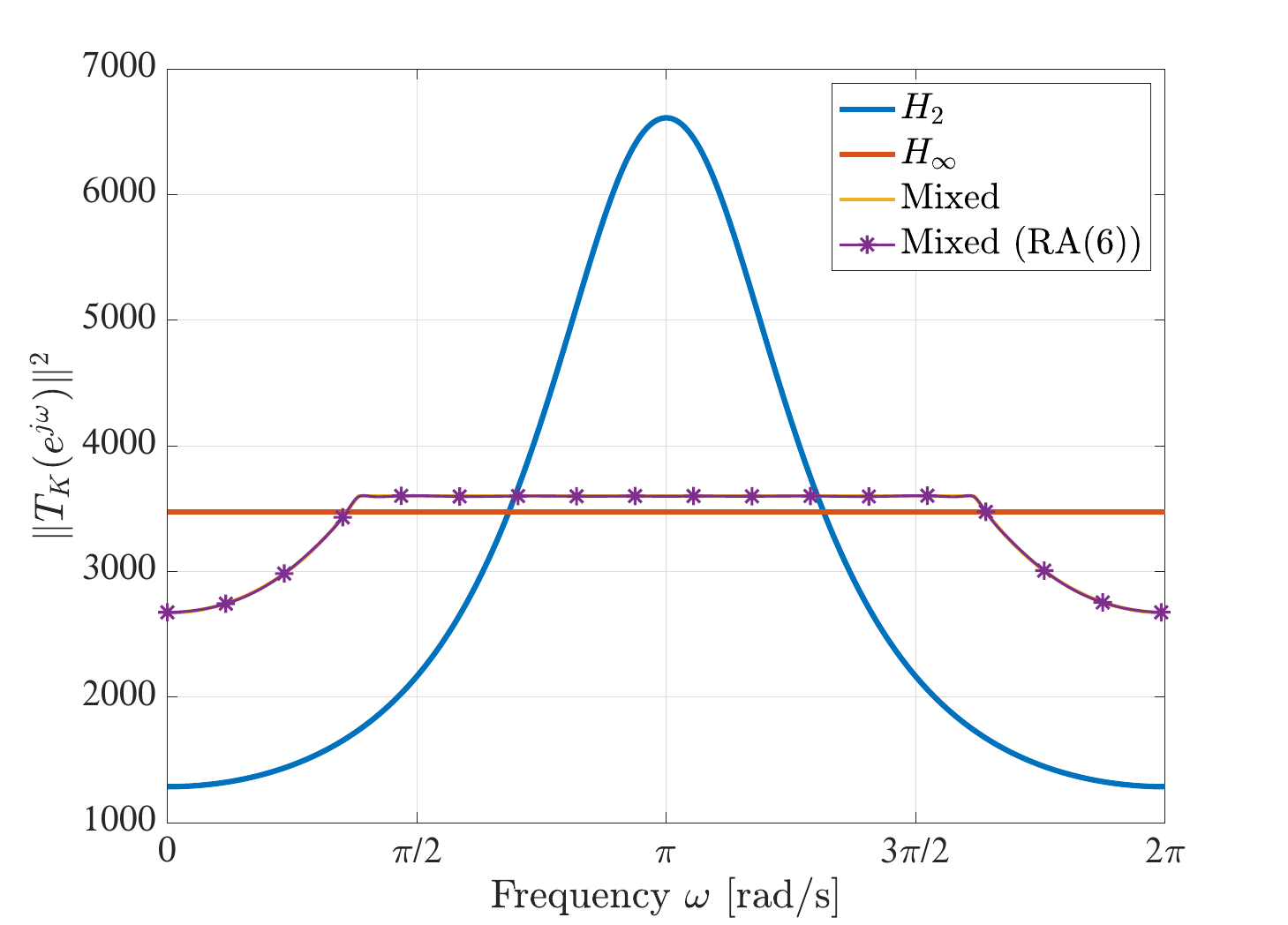}

	\caption{The spectral norm, $\overline{\sigma}(T_K^\ast(\ejw)T_K(\ejw))$ of the mixed $\Htwo/\Hinf$ controller ($\gamma = 60$) and a $6^{th}$ order rational approximation at different frequency values, for the system [AC17]. The cost of the rational controller closely follows the optimal mixed $\Htwo/\Hinf$ controller.} 
 
	\label{fig:ac17_approx}
\end{figure}

In this section, we analyze the performance of the mixed $\Htwo/\Hinf$ controller and the rational approximation method for 2 systems. We use benchmark models from \cite{aircraft}. In particular, we test the aircraft system [AC17] and a chemical reactor system [REA4]. First, we present frequency domain plots for each system, highlighting their performances. The plots show the performance of the mixed $\Htwo/\Hinf$ and the rational approximations. Additional data is presented in Tables \ref{table:ac17} and \ref{table:rea4}. Next, we provide numerical evidence supporting the exponential convergence of the proposed algorithm Algorithm \ref{alg:fixed_point}. Note that in this section, for brevity, a rational approximation of order $m$ is denoted $RA(m)$.


We first present the performance metrics for the [AC17] system. For this system, $\ghi = 58.94$ and $\ght = 81.309$. Thus, the $\Hinf$ norm of the mixed $\Htwo/\Hinf$ controller $\gamma \in (58.94,81.309]$. The system is a $4^{th}$ order system. The system matrices are,
\begin{equation}
A=\left[\begin{array}{cccc}
-2.98 & .93 & 0  & -.034 \\
-.99 & -.21  & .035 & -.001 \\
0 &0 &0 & 1 \\
.39 & -5.55 & 0 & -1.89 
\end{array}\right]  B_u=\left[\begin{array}{cc}
-.032 \\
0 \\
0 \\
-1.6
\end{array}\right] \nonumber.
\end{equation}
We present the frequency domain plots in Figures \ref{fig:ac17} and \ref{fig:ac17_approx}. Note that $\|\T_\K \|_{\op}^2 = \max_{0 \leq \omega \leq 2\pi} \overline{\sigma}( T_K^\ast(\ejw)T_K(\ejw))$. This metric is visualised in Figure \ref{fig:ac17} which highlights the mixed nature of the mixed $\Htwo/\Hinf$ controller. For $\gamma$ close to $\ghi$, the mixed $\Htwo/\Hinf$ controller closely follows the $H_{\infty}$ controller for most of the frequencies but still has less area under the curve ($\norm[\HS]{\T_{\K}}$) and seems to follow the $H_2$ controller. As we increase $\gamma$, we see that the mixed $\Htwo/\Hinf$ controller clips the value of  $\|\T_\K \|_{\op}$ at $\gamma$ for some frequencies, and for the other frequencies, it follows the $H_2$ controller which is to minimise $\norm[\HS]{\T_{\K}}$.
We now focus on the performance of the rational approximations of the mixed $\Htwo/\Hinf$ controller. Since we approximate the spectrum $\NN$, the order of controller is given by the order of the system plus the order of the rational approximation. Table \ref{table:ac17} highlights the performance metrics of the rational approximations. Moreover, Figure \ref{fig:ac17_approx} showcases how a rational approximation looks in the frequency domain. One can observe from Figure \ref{fig:ac17_approx} that $\|\T_\K \|_{\op}$ for the rational approximation can be slightly higher than that of the actual mixed $\Htwo/\Hinf$ controller. As can be seen from Table \ref{table:ac17}, a higher order approximation results in a controller with performance metrics close the the optimal mixed $\Htwo/\Hinf$ controller. For the [AC17] system, a $6^{th}$ order approximation of the spectrum $\NN$ provides seems to be a good choice for the controller.

\begin{table}[htbp]
    \centering
    \setlength\tabcolsep{5pt} 
    \begin{tabular}{|c||c|c||c|c|} 
    \hline
    \textbf{ } & \multicolumn{2}{c||}{$\gamma = 60$}  & \multicolumn{2}{c|}{$\gamma = 68$} \\
        \hline
        \textbf{ } & \textbf{$\norm[\HS]{\T_{\K}}$} & \textbf{$\norm[\op]{\T_{\K}}$} & \textbf{$\norm[\HS]{\T_{\K}}$} & \textbf{$\norm[\op]{\T_{\K}}$} \\
        \hline \hline
        \textbf{$H_{\infty}$} & 58.94 & 58.94 & 58.94 & 58.94  \\
        \hline
        \textbf{Mixed $\Htwo/\Hinf$} & 57.92 & 60 & 54.94 & 68 \\
        \hline
        \textbf{RA(1)} & 58.14 & 60.36 & 54.94 & 69.46 \\
        \hline
        \textbf{RA(3)} & 58.04 & 60.42 & 54.95 & 68.31  \\
        \hline
        \textbf{RA(6)} &  57.92 & 60.07 & 54.94 & 68.07 \\
        \hline
        \textbf{$H_2$} & 54.28 & 81.309 & 54.28 & 81.309  \\
        \hline
    \end{tabular}
        \caption{The performance characteristics of the mixed $\Htwo/\Hinf$ controller obtained from degree 1, 2, and 3 rational approximations to $N(e^{j\omega})$.}
    \label{table:ac17}
\end{table}

\begin{figure}[t]
	\centering
		\centering
		\includegraphics[width=0.4\textwidth]{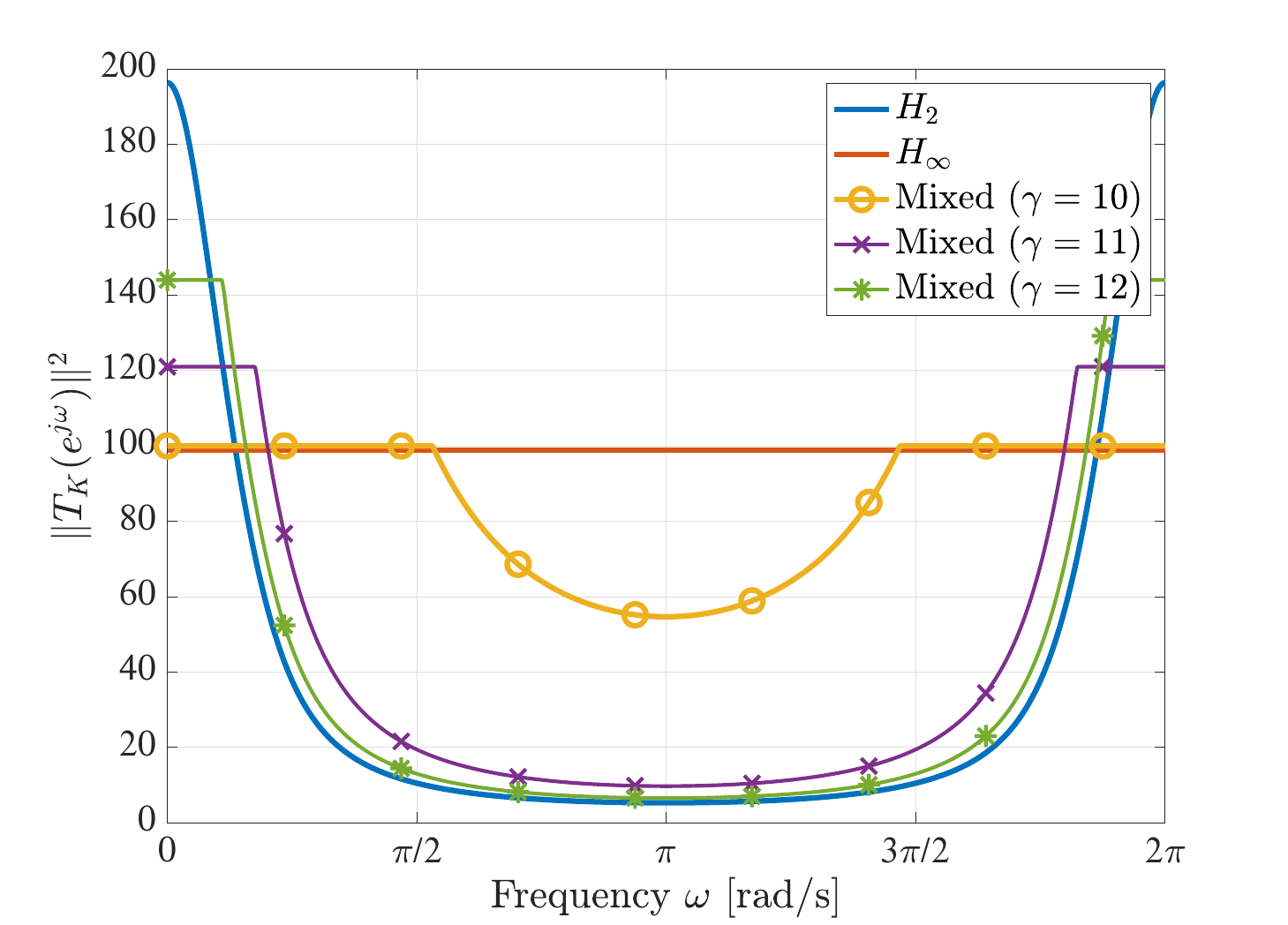}

	\caption{The spectral norm, $\overline{\sigma}(T_K^\ast(\ejw)T_K(\ejw))$ of the mixed $\Htwo/\Hinf$ controller for $\gamma \in \{10, 11, 12\}$ at different frequency values, for system [REA4]. The cost of the mixed $\Htwo/\Hinf$ controller follows $H_2$ and clips at the threshold $\gamma$ for some frequencies.} 
 
	\label{fig:rea4}
\end{figure}

\begin{figure}[t]
	\centering
		\centering
		\includegraphics[width=0.4\textwidth]{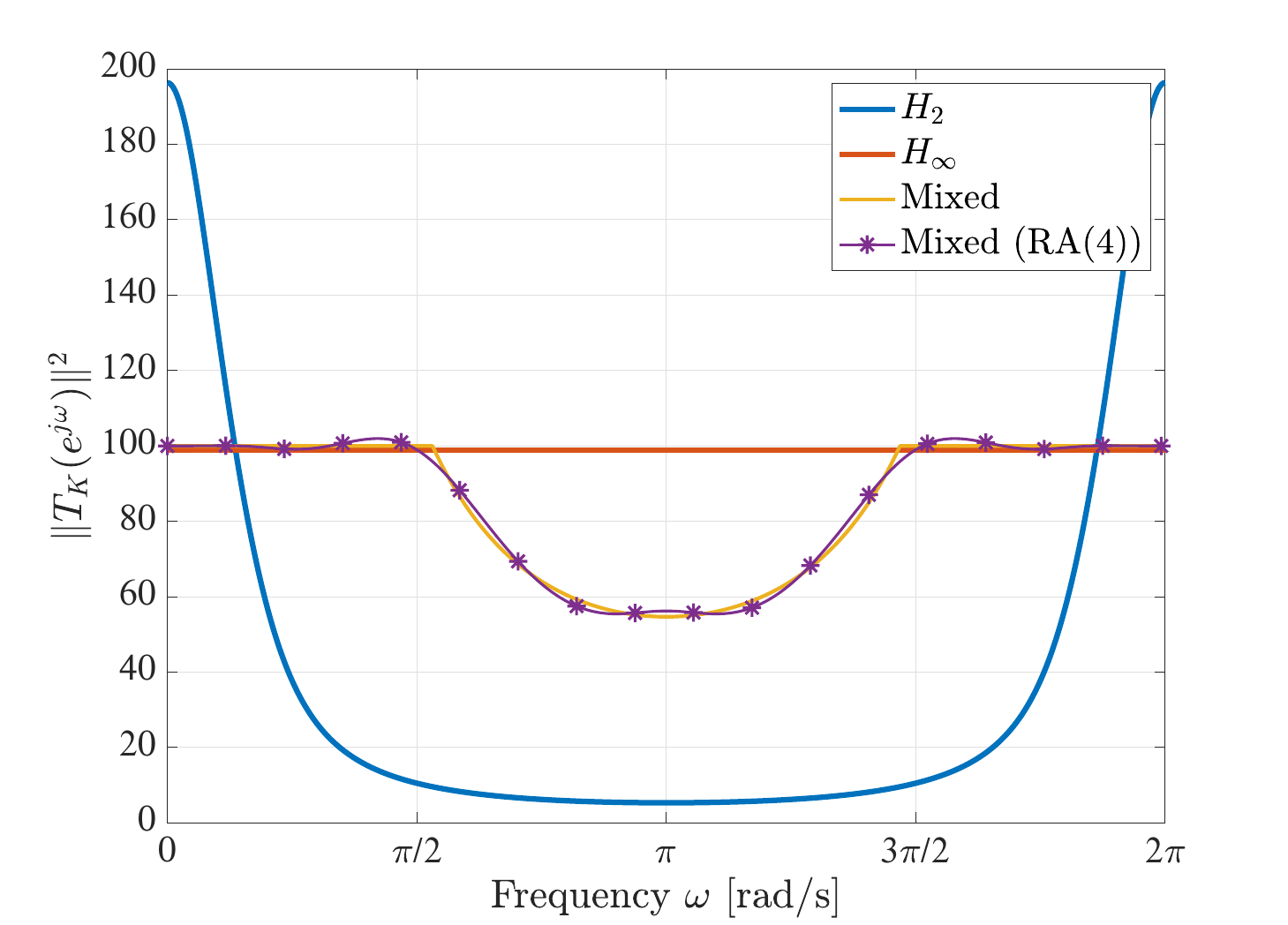}

	\caption{The spectral norm, $\overline{\sigma}(T_K^\ast(\ejw)T_K(\ejw))$ of the mixed $\Htwo/\Hinf$ controller ($\gamma = 10$) and a $4^{th}$ order rational approximation at different frequency values, for system [REA4]. The cost of the rational controller closely follows the optimal mixed $\Htwo/\Hinf$ controller.} 
 
	\label{fig:rea4_approx}
\end{figure}

\begin{figure}[t]
	\centering
		\centering
  \includegraphics[width=0.4\textwidth]{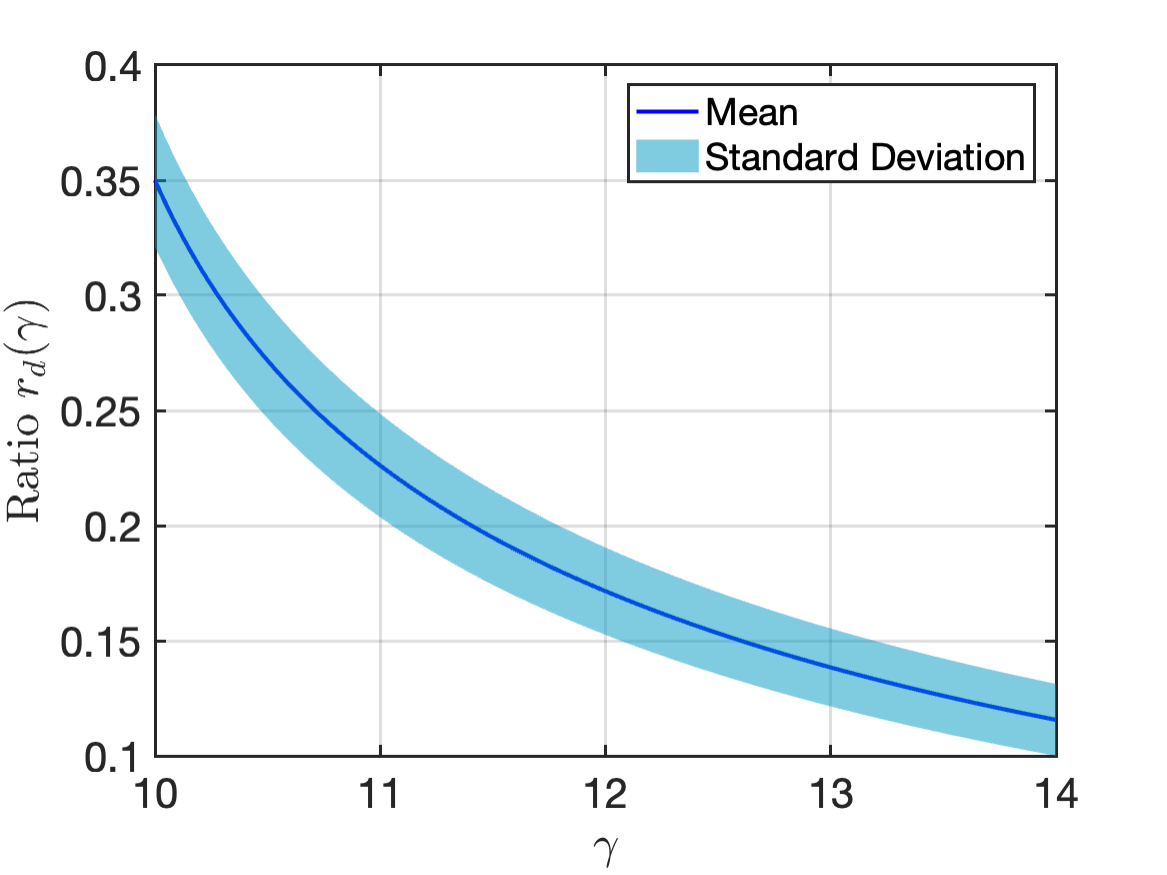}
	\caption{The variation of $r_d(\gamma)$ (defined in \eqref{eq:ratio}) with $\gamma$ for the [REA4] system. The plot indicates different $\NN_1, \NN_2$ in \eqref{eq:ratio} chosen at random. Note that the contraction ratio is always less than $1$ and decreases with an increase in $\gamma$.} 

	\label{fig:contraction}
\end{figure}

We now present the numerical results for the system [REA4]. This is an $8^{th}$ order system which is a chemical reactor system. The system matrices can be found in \cite{aircraft}. As in the case of [AC17], we see in Figure \ref{fig:rea4}, that the mixed $\Htwo/\Hinf$ controller closely follows the $H_{\infty}$ controller for most of the frequencies but still has less area under the curve ($\norm[\HS]{\T_{\K}}$) and seems to follow the $H_2$ controller for some part of the frequencies. The performance of a controller obtained via a $4^{th}$ order rational approximation of the spectrum $\NN$ is shown in Figure \ref{fig:rea4_approx}. As expected, a higher order approximation results in a controller that better approximates the optimal mixed $\Htwo/\Hinf$ controller. For the [REA4] system, a controller obtained via an $8^{th}$ order approximation well approximates the optimal mixed $\Htwo/\Hinf$ controller, as shown in Table \ref{table:rea4}.

\begin{table}[htbp]
    \centering
    \setlength\tabcolsep{5pt} 
    \begin{tabular}{|c||c|c||c|c|} 
    \hline
    \textbf{ } & \multicolumn{2}{c||}{$\gamma = 10$}  & \multicolumn{2}{c|}{$\gamma = 12$} \\
        \hline
        \textbf{ } & \textbf{$\norm[\HS]{\T_{\K}}$} & \textbf{$\norm[\op]{\T_{\K}}$} & \textbf{$\norm[\HS]{\T_{\K}}$} & \textbf{$\norm[\op]{\T_{\K}}$} \\
        \hline \hline
        \textbf{$H_{\infty}$} & 9.94 & 9.94 & 9.94 & 9.94  \\
        \hline
        \textbf{Mixed $\Htwo/\Hinf$} & 9.2 & 10 & 6.17 & 12 \\
        \hline
        \textbf{RA(2)}& 9.3 & 10.05 & 6.21 & 12.09 \\
        \hline
        \textbf{RA(4)} & 9.21 & 10.1 & 6.177 & 12.02  \\
        \hline
        \textbf{RA(8)} & 9.2 & 10.008 & 6.17 & 12.01 \\
        \hline
        \textbf{$H_2$} & 6.06 & 14.01 & 6.06 & 14.01  \\
        \hline
    \end{tabular}
        \caption{The performance characteristics of the mixed $\Htwo/\Hinf$ controller obtained from degree 1, 2, and 3 rational approximations to $N(e^{j\omega})$.}
    \label{table:rea4}
\end{table}

We now present numerical evidence that suggests an exponential convergence of our iterative method Algorithm \ref{alg:fixed_point}. In Figure \ref{fig:contraction}, we consider the ratio,
\begin{align}
    r_d(\gamma) = \frac{d(\NN^{(1)}_1, \NN^{(1)}_2)}{d(\NN_1^{(0)}, \NN_2^{(0)})}.
    \label{eq:ratio}
\end{align}
Here, $\NN_1(e^{j\omega})$ and $\NN_2(e^{j\omega})$ are the frequency domain representation of the power spectra $\NN_1$ and $\NN_2$. $N^{(1)}(e^{j\omega})$ is the spectrum obtained after one iteration of Algorithm \ref{alg:fixed_point} initialized with $N(e^{j\omega})$, and $d(\cdot)$ is a distance metric. For our results, we consider the $H_{\infty}$ norm given by,
\begin{align}
    d(\NN_1, \NN_2) = \max_{\omega \in [0, 2\pi]} \overline{\sigma}(N_1(e^{j\omega}) - N_2(e^{j\omega})).
\end{align}

We consider various random initializations of the $\NN_1$ and $\NN_2$. What we would like to observe is $r_d(\gamma) < 1$ for all $\gamma > \ghi$. 
Since this would mean that after each iteration, the spectrum is close to the optimal solution. As can be seen in Figure \ref{fig:contraction}, the contraction ratio is indeed a decreasing function of $\gamma$ and is always less than $1$. This suggests that Algorithm \ref{alg:fixed_point} is exponentially convergent.

\section{Conclusion}
In this paper, we studied the problem of mixed $\Htwo/\Hinf$ control in the infinite-horizon setting. We provide the exact closed-form solution to the infinite-horizon mixed $\Htwo/\Hinf$ control in the frequency domain. Despite being non-rational, we show that the optimal controller admits a finite-dimensional parameterization. Leveraging this fact, we introduce an efficient iterative algorithm that finds the optimal causal controller in the frequency domain. We show that this algorithm is convergent when the system is scalar and present numerical evidence for exponential convergence of the proposed algorithm. To obtain a finite order controller, we use a rational approximation method (based on the $H_\infty$ norm) and present its performance. In future works, we will extend our results to cases when $d_w > 1$. As mentioned in the paper, numerical results hint that the algorithm is exponentially convergent. Future works will involve analyzing the convergence properties of the proposed iterative algorithm.

\bibliographystyle{IEEEtran}
\bibliography{refs}





\end{document}